# Geometry-based Multi-beam Survey Line Layout Problem

Chuangqi Li[&*ac], Yuhang Wang[&ad], Fan Hu[&b]

[a]College of Computer and Information Engineering, Henan Normal University

[b]College of Physics, Henan Normal University

[c]Key Laboratory of Artificial Intelligence and Personalized Learning in Education of Henan Province

[d]Engineering Lab of Intelligence Business and Internet of Things of Henan Province

[*] Corresponding author: lichuangqi@stu.htu.edu.cn

[&] These authors contributed equally to this work and should be considered co-first authors

## Abstract

The multi-beam measurement system plays a crucial role in ocean mapping and underwater terrain detection. By simultaneously transmitting multiple beams, the system can accurately receive sound waves reflected from the seabed, providing more precise and comprehensive water depth information while effectively revealing the complexity and characteristics of underwater terrain. Building upon the background and application provided by Question B of the 2023 National Mathematical Contest in Modeling for College Students, this paper investigates the relationship between ocean floor width measurement and factors such as beam position, angle, and slope. Utilizing geometric relations, trigonometric similarity, and sine theorem, a mathematical model is established to determine adjacent strip coverage width and overlap ratio. Furthermore, an optimal strategy is determined using a greedy algorithm, while binary search backtracking is employed to derive the interval of the next adjacent survey line with the required overlap ratio in order to obtain an optimal terrain detection strategy.

## 1   Introduction

The multi-beam measurement system represents a significant technological innovation in recent years within the field of marine surveying, bearing revolutionary significance. In various domains such as marine geological exploration, underwater pipeline deployment, and marine resource surveying, the precise acquisition of seabed topographic information is of paramount importance.

### 1.1 Background

The multi-beam depth measurement technique, as a method based on sonar for measuring seabed depth, faces critical challenges when dealing with complex marine environments, such as uneven seabed. A key issue in this regard is the design of a rational measurement route that ensures comprehensive coverage of the marine area, captures sufficient data points, and simultaneously maintains a specified overlap ratio between adjacent measurement lines. This poses a central challenge in research. In-depth investigation of this problem through the establishment of corresponding mathematical models contributes to optimizing the implementation of multi-beam depth measurement, enhancing measurement efficiency, and demonstrating practical applicability.

**1.2 Comprehensive Analysis**

The irregularity of seafloor topography presents a significant challenge to the field of Marine surveying. In the depths of the ocean, the terrain can assume various intricate forms, such as hills, ravines, and canyons. These wave features have the potential to obstruct or refract sound waves during their propagation, leading to issues with beam overlap and leakage in measurements. The complexity of seabed topography often has a substantial impact on measurement accuracy and data integrity. To tackle these challenges, this study aims to establish a coordinate system tailored for a specific ocean area and seabed topographic conditions while leveraging geometric relations like triangle similarity and sine theorem[1]. A mathematical model is formulated to determine water depth and coverage width accurately while calculating the overlap ratio of survey ship positions[2]. For addressing the problem of determining the shortest distance between measurement lines, we first prove optimality using a method before employing a greedy algorithm for its solution. Finally, we compare our results with standard answers to evaluate their accuracy. The overall solution framework for this study is illustrated in Figure 1.

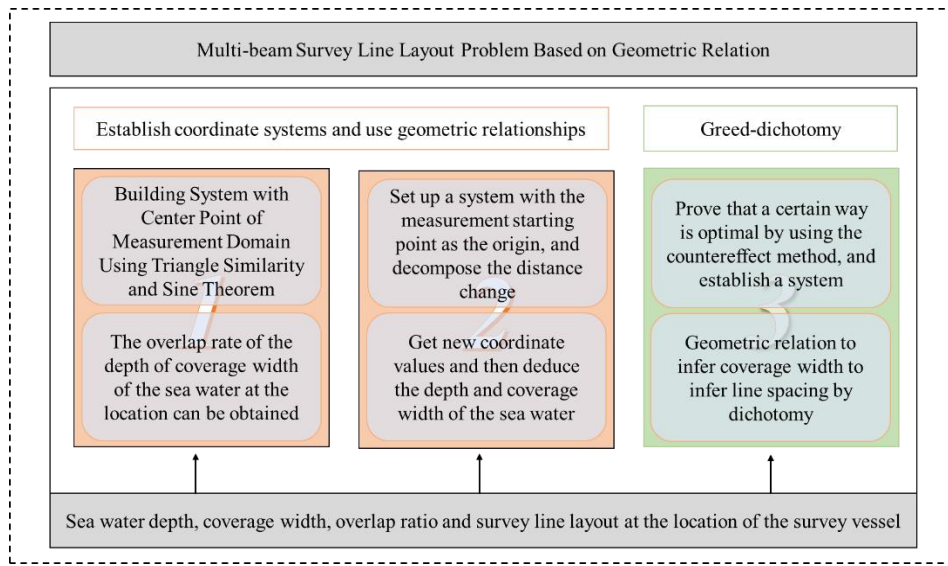

Figure 1: Overall Flowchart

## 2 Model Establishment and Solution

**2.1 Model Establishment**

Given known parameters such as angles and lengths, in order to establish a mathematical model for the coverage width and overlap ratio between adjacent strips in multi-beam depth measurement, geometric relationships, and principles such as the sine theorem and the theorem of triangular similarity are employed. The geometric model is established, and the model results are obtained through C++ programming.

As depicted in Figure 2, the following coordinate system is established: taking a point on the slope surface (as shown in Figure 2) as the origin, considering the projection of the slope surface normal onto the horizontal plane as the x-axis, and establishing a rectangular coordinate system with a y-axis perpendicular to the x-axis. The angle between the survey line direction and the x-axis is denoted as $\beta$.

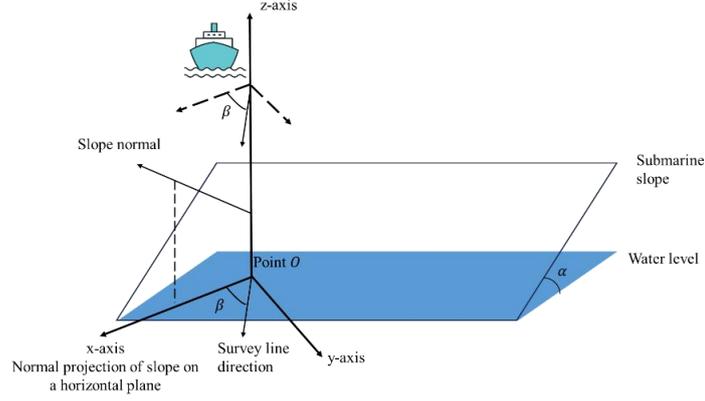

Figure 2: Schematic Diagram of Coordinate Axes for Problem Two

## Step1: The depth of seawater at a specific location

When the distance of the measurement vessel from the center point of the marine area (0, 0, 120) is set as $dist$, its components along the x-axis and y-axis are

$$\begin{cases} dist_x = dist \cdot \cos\beta \\ dist_y = dist \cdot \sin\beta \end{cases} \quad (1)$$

Hence, the coordinates of this point are $(dist_x, dist_y, 120)$. Additionally, since the x-axis is parallel to the horizontal plane, the specific relationship is illustrated in Figure 3 as depicted below.

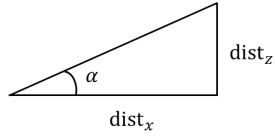

Figure 3: Abstract Two-Dimensional Diagram with Respect to the Horizontal Plane

So, the variation of the measurement vessel in the z-direction is

$$dist_z = \tan\alpha \cdot dist_x \quad (2)$$

Therefore, the seawater depth at this location should be

$$D = 120 + dist_z \quad (3)$$

## Step2: The new slope angle $\gamma$ at a certain location

The directional vector of the survey line is $\vec{n_1} = (\cos\beta, \sin\beta, 0)$, and the coordinates of the measurement vessel at its arrival position are $A(dist_x, dist_y, 120)$. The equation of the plane with the survey line direction as the normal vector passing through the point where the measurement vessel is located is given by

$$\begin{cases} \vec{n_1} \cdot \overrightarrow{AB} = 0 \\ (\cos\beta, \sin\beta, 0) \cdot (x - dist_x, y - dist_y, z - dist_z) \end{cases} \quad (4)$$

and $B(x, y, z)$ represents a point on this plane. This study deduced through simplification

$$x\cos\beta + y\sin\beta - (\cos\beta\, dist_x + \sin\beta\, dist_y) = 0 \quad (5)$$

Now, to find the plane equation for the slope surface, let a certain point on the slope surface be the origin, and the normal vector of the slope surface be $\vec{n_2} = (\sin\alpha, 0, \cos\alpha)$. The plane equation for the slope surface is given by formula (6), where $(x, y, z)$ is a point on the slope surface.

$$\begin{cases} \vec{n_2} \cdot (x, y, z) = 0 \\ x\sin\alpha + z\cos\alpha = 0 \end{cases} \quad (6)$$

Through the above reasoning, this study derived the equations of the two planes above. According to mathematical principles, the direction vector of the intersection line of two planes is equal to the cross product of their normal vectors, which can be expressed as:

$$\vec{n_3} = \vec{n_1} \times \vec{n_2} = (\sin\beta\cos\alpha, -\cos\beta\cos\alpha, -\sin\beta\sin\alpha) \quad (7)$$

The directional vector of the intersection line's projection onto the horizontal plane is $\vec{n_4} = (\sin\beta\cos\alpha, -\cos\beta\cos\alpha, 0)$. It is evident that the cosine value of the angle formed by the intersection line and the bottom surface is

$$\cos\gamma = \frac{\vec{n_3} \cdot \vec{n_4}}{|\vec{n_3}| \cdot |\vec{n_4}|} = \frac{\cos\alpha}{\sqrt{\cos^2\alpha + \sin^2\beta \sin^2\alpha}} \quad (8)$$

Furthermore, since the angle $\gamma$ must be less than 90°, and $\gamma$ represents the new slope angle, and the seabed depth at the required location has been determined in the previous context, thus, the coverage width can be calculated based on the following formula:

$$W = 120 \cdot \sin\frac{\theta}{2} \cdot \left[ \frac{1}{\sin\left(\frac{\pi}{2} - \frac{\theta}{2} - \gamma\right)} + \frac{1}{\sin\left(\frac{\pi}{2} - \frac{\theta}{2} + \gamma\right)} \right] \quad (9)$$

## 2.2 Model Solution

The solution of the model involves simultaneously solving the aforementioned equations by substituting the known parameters directly. A step-by-step process leads to the results listed in Table 1.

Table 1: Relationship Table between Survey Line Angle and Coverage Width

| Coverage Width/m | Distance of the Measurement Vessel from the Center Point of the Marine Area / Nautical Miles | | | | | | | |
| --- | --- | --- | --- | --- | --- | --- | --- | --- |
| | 0 | 0.3 | 0.6 | 0.9 | 1.2 | 1.5 | 1.8 | 2.1 |

| Survey Line Direction Angle/° | 0 | 415.69 | 466.09 | 516.49 | 566.88 | 617.28 | 667.68 | 718.08 | 768.48 |
|---|---|---|---|---|---|---|---|---|---|
| | 45 | 415.69 | 451.33 | 486.96 | 522.60 | 558.24 | 593.87 | 629.51 | 665.15 |
| | 90 | 415.69 | 415.69 | 415.69 | 415.69 | 415.69 | 415.69 | 415.69 | 415.69 |
| | 135 | 415.69 | 380.05 | 344.41 | 308.78 | 273.14 | 237.50 | 201.86 | 166.23 |
| | 180 | 415.69 | 365.29 | 314.89 | 264.49 | 214.09 | 163.69 | 113.29 | 62.90 |
| | 225 | 415.69 | 380.05 | 344.41 | 308.78 | 273.14 | 237.50 | 201.86 | 166.23 |
| | 270 | 415.69 | 415.69 | 415.69 | 415.69 | 415.69 | 415.69 | 415.69 | 415.69 |
| | 315 | 415.69 | 451.33 | 486.96 | 522.60 | 558.24 | 593.87 | 629.51 | 665.15 |

As the distance of the measurement vessel from the center point of the marine area decreases, the variation in coverage width becomes smaller for different survey line direction angles[3]. This aligns with intuitive understanding and, to a certain extent, indicates the correctness of the model. Due to the existence of symmetry, if the sum of two survey line direction angles is 360°, the values of coverage width should be equal.

In a practical application scenario, consider a rectangular marine area with a length of 2 nautical miles in the north-south direction and a width of 4 nautical miles in the east-west direction. The seawater depth at the center point of the marine area is 110 m, with a slope inclination of 1.5° (deeper in the west and shallower in the east). A multi-beam transducer with an opening angle of 120° can completely cover the entire surveyed area, and the overlap ratio between adjacent strips meets the requirement of 10% to 20%.

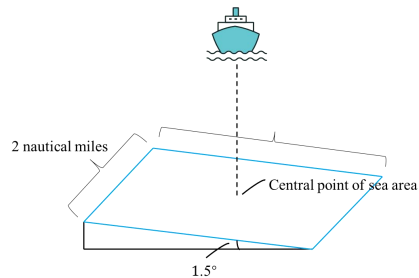

Figure 4: Application Topographic Map of the Real-world Problem

In this study, the greedy strategy is that the measurement lines are parallel in the north-south direction, from deep water to shallow water, that is, from west to east. Each measurement line is positioned exactly 10 percent over the previous line. This process continues until the entire survey area is covered[4][5].

As illustrated in Figure 5, establish the following coordinate system with the eastward direction as the positive axis. Given that the seawater depth at the center point of the marine area is 110m, the length of segment $IZ$ is 110m. Let the length of segment $ZP$ be $D_1$, and designate the seawater depth at the point where the measurement vessel is located as $D$. The perpendicular line from the measurement vessel to $PT$, intersecting at point $Q$ [6].

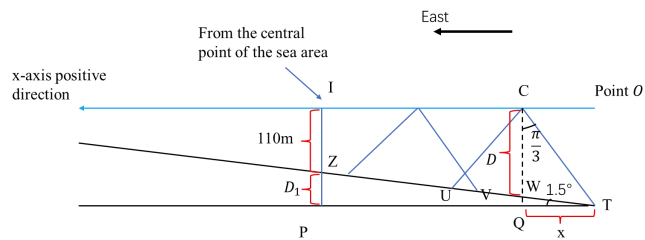

Figure 5: Schematic Diagram of the Cross-Sectional Profile of the Marine Area

**Step1: Determining the Starting Point and Locating.**

the Seawater Depth at the Measurement Vessel Location

As the point $I$ serves as the center of the marine area, with the map oriented in the east-west direction, the length of segment $PT$ is determined to be 3704 meters. Considering the slope inclination of 1.5°, a compensatory adjustment is applied between $\triangle PTZ$ and $\triangle WTQ$, leading to the determination of the seawater depth at the measurement vessel location[7]. The seawater depth at the location of the measurement vessel is thus

$$D = 110 + D_1 - x \tan \alpha \tag{11}$$

**Step2: Calculating Coverage Width from the Seawater Depth at the Given Point**

In the context of $\triangle CWT$ and $\triangle CWU$, applying the Law of Sines allows us to determine $\triangle CWU$ and $|WU|$. Therefore, the coverage width is

$$W = |WT| + |WU| = D \times \left( \frac{\sin 60°}{\sin 28.5°} + \frac{\sin 60°}{\sin 31.5°} \right) \tag{12}$$

**Step3: Using Binary Search to Infer the Interval Between Survey Lines Based on the Overlap Ratio**

The formula for the overlap ratio, as depicted in equation (13), allows us to derive the expression for d.

$$d = (1-\eta) \times \overline{w}, \overline{w} = (\frac{w_1 + w_2}{2}) \tag{13}$$

The problem assumes an overlap ratio of $\eta = 10\%$. By employing $\eta$ with a binary search, the position of the next survey line is inversely determined. The binary search method is utilized for its advantages in speed, simplicity, and ensured convergence, requiring low demands on the function's characteristics, making it highly implementable. As the optimal strategy has been previously demonstrated, the application of a greedy algorithm for direct resolution is viable. The algorithm involves computing the initial coordinates of the selected survey lines, entering a loop for iteration[8][9][10]. If the condition is met, the iteration concludes, outputting the number of survey lines. Otherwise, the loop continues.

Upon obtaining the formula for the interval between adjacent survey lines, continuous iterations enable the solution to be derived as follows: $D_1 = 96.9927$, The coordinates of different survey line positions, the overlap ratio between the current survey line and the preceding one, and the coverage width at the current position are presented in the following table 2:

Table 2: Results for Problem Three

| The current survey line position x | the overlap ratio with the preceding survey line | The coverage width at the current position. |
|---|---|---|
| 358.522 | / | 686.168 |
| 951.734 | 0.10001 | 632.228 |
| 1498.31 | 0.10001 | 582.528 |
| ⋮ | ⋮ | ⋮ |
| 7308.43 | 0.10001 | 54.2213 |
| 7355.3 | 0.10001 | 49.9589 |
| 7398.49 | 0.10001 | 46.0316 |

Namely, 34 survey lines, totaling 68 nautical miles.

From Table 2, it is observed that as $x$ increases, corresponding to the direction towards shallower waters, the coverage width should decrease. This aligns with the computed results in Table 2, demonstrating the reliability of the outcomes and the accuracy of the model. With the increase of $x$, the distance between adjacent survey lines diminishes. Utilizing the coordinate values from the above table, a coordinate system can be established, and the resulting graph is illustrated in Figure 7.

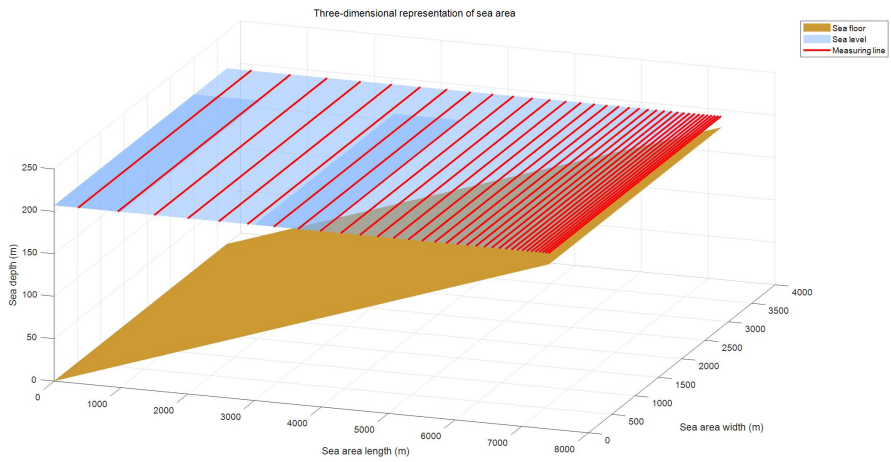

Figure 6: 3D Schematic Diagram of Survey Lines

In the illustration, the yellow background represents the seabed plane, the light blue horizontal plane signifies the sea level, and the red lines depict the direction of the survey lines. This three-dimensional visualization enables a clearer understanding of the seabed, sea level, and the configured survey lines, providing a more intuitive representation of the survey line layout.

## 3 Peroration

In this paper, this study successfully applies the binary algorithm and greedy algorithm to solve the optimization problem of a multi-beam measurement system. The dichotomous algorithm demonstrates good performance in reverting the interval of the next adjacent line with a qualified overlap rate. By gradually reducing the search space, the optimal solution is quickly found, thereby improving the speed and accuracy of problem-solving. This is crucial for determining the appropriate measurement line spacing to ensure that the coverage width and overlap rate of adjacent strips are within the ideal range. The greedy algorithm plays a significant role in determining the shortest line distance. By selecting the current seemingly optimal solution at each stage, the greedy algorithm can achieve the global optimal solution on the whole, providing strong support for the optimization strategy. In a multi-beam measurement system, this helps to effectively arrange the layout of the measurement lines, resulting in more uniform and comprehensive coverage of the entire sea area. Although the dichotomous algorithm and greedy algorithm have clear advantages in solving problems, the implementation of the dichotomous algorithm may be impacted by uneven data distribution and measurement errors, making it difficult to achieve optimal results in certain scenarios.

Drawing on information provided by Problem B of the 2023 National Mathematical Contest in Modeling for College Students, this article uncovers the intricate relationship between seafloor width measurements and the position, angle, and slope of the emitted beam. this article's geometry-based mathematical model, which combines trigonometric similarity and sine theorem, not only provides valuable insights into adjacent strip coverage widths and overlap ratios but is also rigorously tested against real-world data to ensure reliability. This integrated approach supports the implementation of the greedy algorithm and binary search technology to determine the optimal strategy, enhancing the speed and accuracy of problem-solving and ensuring uniformity in measurement line layout.

Furthermore, this study will continue to explore ways to deepen and combine other algorithms to better adapt to diverse terrain and data conditions. The efficiency and applicability of our model are demonstrated through its contribution to improving the practicality and effectiveness of multi-beam measurement systems. The comprehensive research presented in this paper advances ocean mapping techniques while advancing discussions on sustainable ocean exploration and resource management.